\documentclass{article}
\usepackage{graphicx} 
\usepackage{amsmath} 
\usepackage{amssymb}
\title{Short proof of Halin's grid theorem for thick ends}
\author{Ye Chern\thanks{1254015070@ybu.edu.cn}}
\date{September 2025}

\begin{document}

\maketitle
\begin{abstract}
\hspace{2em}In this paper, we present a short proof of Halin's grid theorem.
\end{abstract}
\section{Introduction}
\hspace{2em}Halin's theorem characterizes thick ends. It states that an infinite graph G has a thick end if and only if the hexagonal quarter grid $H^{\infty}$\footnote{We followed R.Diestel,see [1] chapter 8 section 2} is a topological minor of $G$. In our proof, we first verify a special case of this theorem, then deduce the general case.
\section{Proof}
\hspace{2em}Let G be an infinite graph.\\
\textbf{Lemma}: Suppose that there are disjoint rays $C, L_0, L_1,\dotsc$ in graph $G$ such that for each j there are infinitely many $L_j-C$ paths, denoted $P_j$, and $\forall m,n\in \mathbb{N}$, $L_n$ intersects only finitely many paths in $P_m$, then $H^{\infty}$ is a topological minor of $G$.\\
\textit{Proof}: We construct by induction a sequence of graphs $H_n$ and finite subsets $C_n\subset C$. Initially, recall that for any ray there is a natural order on it, let $H_0$ be the union of $L_0$, C, and a path $l$ in $P_0$ such that $l\cap L_0$ has the minimum element in $C\cap P_0$, $C_0$ be that vertex.\\
In the induction step, suppose that we already construct $H_n$, let us construct $H_{n+1}$. We choose a sub-ray L of $L_{n+1}$ such that it is disjoint with all the paths in $P_j$ for $j\leq n$. Since there are infinitely many $L_j-C$ paths for each j, we can choose a point $x_0\in C$ such that it is strictly greater than any point in $C_n$ and there is a $L_0-C$ path $l_0=a_0\dots x_0$.Suppose we have chosen $x_j$ and $l_j$ for $j<k\leq n+1$, we choose $x_k$ such that it is strictly greater than any point in $C_n\cup \{x_0, x_1,\dots, x_{k-1}\}$ and there is a $L_k-C$ path $l_k=a_k\dots x_k$.We define $H_{n+1}=H_n\cup L_{n+1}\cup l_0\cup l_1\dots \cup l_{n+1}$, and $C_{n+1}=C_n\cup \{x_0, x_1,\dots ,x_{n+1}\}$, which is finite.\\
Let H=$\lim_{n\to \infty}H_n=\cup _{n=0}^{\infty}H_n\subset G$ and $C'=\cup _{n=0}^{\infty}C_n$. It is easy to see that H has $H^{\infty}$ as a topological minor. We only need to choose all the points at odd indices in $C'\cap L_0$(with respect to the natural order on $L_0$), and connect $L_0$ with C via these paths corresponding to these points, and connect with $L_1$ via C, then connect $L_0$ with C via these paths corresponding to these points that has even indices in $C'\cap L_1$, and so on. Since $H\subset G$, we have proved the lemma.$\hfill\blacksquare$\\
\textbf{Halin's grid theorem}: For an infinite graph $G$, $G$ has a thick end iff $H^{\infty}$ is a topological minor of $G$.\\
\textit{Proof}: The 'If' part is trivial. We now prove the ‘only if’ part.Let $\pi$ be a thick end of G.\\
We use proof by contradiction. Suppose this theorem is false, we define by induction a sequence of graph $F_n$ whose vertices are rays in $\pi$ and edges are sets of paths. In the beginning, choose a ray $R_0$ in $\pi$, let $F_0=\{R_0\}$. Suppose that we already construct $F_{n-1}$, let us construct $F_{n}$. Let $F_{n-1}=\{R_0,R_1,\dots ,R_{n-1}\}$, and $\mathcal{E}_{ij}$ be some $R_i-R_j$ paths by our choice if $R_i$ and $R_j$ are adjacent and null set if they are not are adjacent. Since $\pi$ is thick, we can choose $R_n \notin F_{n-1}$. If for any $\mathcal{E}_{ij}$, $R_{n}$ intersects only with finitely many paths in $\mathcal{E}_{ij}$, then we choose an infinite set of disjoint $R_{n-1}-R_n$ path A. If for some j, $R_j$ is intersects with infinitely many elements in A, let A' be such paths' set. For each $R_{n-1}-R_{n}$ path L in A', suppose its last point in $R_j$is y, then yL is $R_j-R_n$ path, and does not intersect with $R_j$ except at y. So we can obtain an infinite set of disjoint $R_j-R_n$ paths. Repeat our process for such set. This process will finish in finitely many steps since $F_n$ has only n elements. And when it stop, suppose we get $R_i-R_n$ paths finally, then we define $F_n$ to be $F_{n-1}\cup \{R_n\}$ with $R_n$ is adjacent to $R_i$ only. If for some $\mathcal{E}_{ij}$, $R_n$ intersects with infinitely many paths in $\mathcal{E}_{ij}$, then we add a vertex in the midpoint of the line segment $\mathcal{E}_{ij}$ for every such $\mathcal{E}_{ij}$, and then identify all the midpoints we add(quotient graph). We define $F_n$ is the result graph.\\
In our process of construct of $F_n$, if a line segment will be add point infinitely, identify the line segment with the interval [0,1], by elementary analysis, we can see there is a sequence of number $a_0<a_1<a_2<\dots $, each number represent a ray, and $a_n$ is adjacent to $a_{n+1}$, so it is easy to see it contains a $TH^{\infty}$ by construction, which is contradiction to our assumption, so it is impossible. Hence, each edge is completely determined after finitely many step. So we can consider the graph $F=\lim_{n \to \infty }F_n$, which is clearly an infinite and connected graph. By the star-comb lemma \footnote{see [1] chapter 8 lemma 8.2.2}, it contains either a star or a comb. The star case is our lemma proved above, and the comb case is obvious since it contains a ray.
$\hfill\blacksquare$
\section{References}
[1] Diestel, R. (2017). \textit{Graph Theory} (5th ed.). Springer
\end{document}